\documentclass[12pt,twoside]{article}

\usepackage{amsfonts}
\usepackage{amssymb}
\usepackage{amsmath}
\usepackage{graphicx}

\newcommand{\ignore}[1]{}

\def\@begintheorem#1#2{\par\bgroup{\sc #1\ #2. }\it\ignorespaces}
\def\@opargbegintheorem#1#2#3{\par\bgroup{\sc #1\ #2\ (#3). } \it\ignorespaces}
\def\@endtheorem{\egroup}
\newtheorem{theorem}{Theorem}[section]
\newtheorem{corollary}[theorem]{Corollary}
\newtheorem{lemma}[theorem]{Lemma}
\newtheorem{example}[theorem]{Example}
\newtheorem{proposition}[theorem]{Proposition}
\newtheorem{definition}[theorem]{Definition}
\newcommand{\bt}[1]{\begin{theorem}\label{#1}}
\newcommand{\bc}[1]{\begin{corollary}\label{#1}}
\newcommand{\bl}[1]{\begin{lemma}\label{#1}}
\newcommand{\be}[1]{\begin{example}\label{#1}}
\newcommand{\bp}[1]{\begin{proposition}\label{#1}}
\newcommand{\ba}[1]{\begin{algorithm}\rm\label{#1}}
\newcommand{\bd}[1]{\begin{definition}\rm\label{#1}}
\newcommand{\bpr}{\noindent {\em Proof. }}
\newcommand{\et}{\end{theorem}}
\newcommand{\ec}{\end{corollary}}
\newcommand{\el}{\end{lemma}}
\newcommand{\ee}{\end{example}}
\newcommand{\ep}{\end{proposition}}
\newcommand{\ed}{\end{definition}}
\newcommand{\epr}{{\ \vbox{\hrule\hbox{%
\vrule height1.3ex\hskip0.8ex\vrule}\hrule}}\\\par}

\def\Z{\mathbb{Z}}

\def \supp {{\rm supp}}

\begin{document}

\title{\bf Huge Unimodular N-Fold Programs}
\author{
Shmuel Onn
\thanks{\small Technion - Israel Institute of Technology, Haifa, Israel.
Email: onn@ie.technion.ac.il}
\and
Pauline Sarrabezolles
\thanks{\small Ecole des Ponts, Paris, and
Technion, Haifa. Email: pauline.sarrabezolles@gmail.com}
}

\date{}

\maketitle

\begin{abstract}
Optimization over $l\times m\times n$ integer $3$-way tables with given line-sums
is NP-hard already for fixed $l=3$, but is polynomial time solvable with both $l,m$ fixed.
In the {\em huge} version of the problem, the variable dimension $n$ is encoded in
{\em binary}, with $t$ {\em layer types}. It was recently shown that the huge problem can
be solved in polynomial time for fixed $t$, and the complexity of the problem for
variable $t$ was raised as an open problem. Here we solve this problem and show that the
huge table problem can be solved in polynomial time even when the number $t$ of types
is {\em variable}. The complexity of the problem over $4$-way tables with variable $t$
remains open. Our treatment goes through the more general class of
{\em huge $n$-fold integer programming problems}. We show that huge integer programs
over $n$-fold products of totally unimodular matrices
can be solved in polynomial time even when the number $t$ of brick types is variable.

\end{abstract}

\section{Introduction}

Consider the following optimization problem over $3$-way tables with given line-sums:
$$\min\{wx\ :\ x\in\Z_+^{l\times  m\times n}\,,\ \sum_i x_{i,j,k}=e_{j,k}
\,,\ \sum_j x_{i,j,k}=f_{i,k}\,,\ \sum_k x_{i,j,k}=g_{i,j}\}\ .$$
It is NP-hard already for $l=3$, see \cite{DO1}. Moreover, {\em every}
bounded integer program can be isomorphically represented in polynomial time
for some $m$ and $n$ as some $3\times m\times n$ table problem, see \cite{DO2}.
However, when both $l,m$ are fixed, it is solvable in polynomial time
\cite{DHOW,HOW,Onn}, and in fact, in time which is cubic in $n$ and linear in the
binary encoding of $w,e,f,g$, see \cite{HOR}. Assume throughout then that $l,m$ are
fixed, and regard each table as a tuple $x=(x^1,\dots,x^n)$ consisting of $n$ many
$l\times m$ {\em layers}. The problem is called {\em huge} if the variable number $n$ of
layers is encoded in {\em binary}. We are then given $t$ {\em types} of layers, where
each type $k$ has its cost matrix $w^k\in\Z^{l\times m}$, column-sums vector
$e^k\in\Z_+^m$, and row-sums vector $f^k\in\Z_+^l$. In addition, we are given positive
integers $n_1,\dots,n_t,n$ with $n_1+\cdots+n_t=n$, all encoded in binary. A feasible
table $x=(x^1,\dots,x^n)$ then must have first $n_1$ layers of type $1$, next $n_2$ layers
of type $2$, and so on, with last $n_t$ layers of type $t$. The special case of $t=1$
is the case of {\em symmetric} tables, where all layers have the same cost, row and column
sums, and the classical (non-huge) table problem occurs as the special case of $t=n$ and
$n_1=\cdots=n_t=1$. Note that for each $k$, the set of possible layers of type $k$ is
$$
\left\{z\in\Z_+^{l\times m}\ :\ \sum_i z_{i,j}=e^k_j
\,,\ \sum_j z_{i,j}=f^k_i\right\}\ ,
$$
and may have cardinality which is exponential in the binary encoding of $e^k,f^k$.
So it is not off hand clear how to even write down a single table, let alone optimize.

\vskip.2cm
The huge table problem was recently considered in \cite{Onn2}, where it was shown,
combining results of \cite{DHOW,HOW,Onn} on Graver bases and results of \cite{ES,GR} on
integer cones, that it can be solved in polynomial time for fixed $t$. The complexity of
the problem for variable $t$ was raised as an open problem. Here we solve this problem and
show that the huge table problem can be solved in polynomial time even when $t$ is variable.

\bt{ThreeWayLineSum}
The huge $3$-way table problem with a variable number $t$ of types can be solved in
time which is polynomial in $t$ and in the binary encoding of $w^k,e^k,f^k,g,n_k$.
\et

It was moreover shown in \cite{Onn2} that the huge $d$-table problem over
$m_1\times\cdots m_{d-1}\times n$ tables with $m_1,\dots m_{d-1}$ fixed and $n$
variable can also be solved in polynomial time for any fixed number $t$ of types.
Interestingly, we do not know whether Theorem \ref{ThreeWayLineSum} could be extended
to this more general situation, and the complexity of the huge $d$-way table problem
with variable $t$ remains open, already for $3\times3\times3\times n$ tables.

\vskip.5cm
Theorem \ref{ThreeWayLineSum} follows from broader results which we proceed to describe.
The class of {\em $n$-fold integer programming} problems is defined as follows.
The {\em $n$-fold product} of an $s\times d$ matrix $A$
is the following $(d+sn)\times(dn)$ matrix, with $I$ the $d\times d$ identity,
$$A^{[n]}\quad:=\quad
\left(
\begin{array}{cccc}
  I    & I    & \cdots & I    \\
  A    & 0      & \cdots & 0      \\
  0      & A    & \cdots & 0      \\
  \vdots & \vdots & \ddots & \vdots \\
  0      & 0      & \cdots & A    \\
\end{array}
\right)\quad .
$$
The classical $n$-fold integer programming problem is then the following:
\begin{equation}\label{classical-n-fold}
\min\left\{wx\ :\ x\in\Z^{dn}\,,\ A^{[n]}x=b\,,\ l\leq x\leq u\right\}\ ,
\end{equation}
where $w\in\Z^{dn}$, $b\in\Z^{d+sn}$, and $l,u\in\Z_{\infty}^{dn}$ with
$\Z_{\infty}:=\Z\uplus\{\pm\infty\}$. For instance, optimization over
multiway tables is an $n$-fold program, as is explained later on.

Our starting point is the following result on classical $n$-fold integer
programming, established in \cite{DHOW,HOW}, building on results of \cite{AT,HS,SS}.
See the monograph \cite{Onn} for a detailed treatment of the theory
and applications of $n$-fold integer programming.

\bp{Classical}
For fixed matrix $A$, the classical $n$-fold integer programming problem
(\ref{classical-n-fold}) can be solved in time polynomial in $n$ and the
binary encoding of $w,l,u,b$.
\ep
This result holds more generally if the identity $I$ in the definition of $A^{[n]}$
is replaced by another fixed matrix $B$. Moreover, recently, in \cite{HOR},
it was shown that the problem can be solved in time which is cubic in $n$ and
linear in the binary encoding of $w,b,l,u$.

\vskip.2cm
The vector ingredients of an $n$-fold integer program are naturally arranged
in {\em bricks}, where $w=(w^1,\dots,w^n)$ with $w^i\in\Z^d$ for $i=1,\dots,n$,
and likewise for $l,u$, and where $b=(b^0,b^1,\dots,b^n)$ with $b^0\in\Z^d$
and $b^i\in\Z^s$ for $i=1,\dots,n$. Call an $n$-fold integer program {\em huge}
if $n$ is encoded in {\em binary}. More precisely, we are now given $t$
{\em types} of bricks, where each type $k=1,\dots,t$ has its cost
$w^k\in\Z^d$, lower and upper bounds $l^k,u^k\in\Z^d$, and
right-hand side $b^k\in\Z^s$. Also given are $b^0\in\Z^d$
and positive integers $n_1,\dots,n_t,n$ with $n_1+\cdots+n_t=n$, all encoded in binary.
A feasible point $x=(x^1,\dots,x^n)$ now must have first $n_1$
bricks of type $1$, next $n_2$ bricks of type $2$, and so on, with last
$n_t$ bricks of type $t$. Classical $n$-fold integer programming occurs as the
special case of $t=n$ and $n_1=\cdots=n_t=1$, and {\em symmetric} $n$-fold
integer programming occurs as the special case of $t=1$. We show the following.

\bt{optimization_theorem}
Let $A$ be a fixed totally unimodular matrix and consider the huge \break $n$-fold
program over $A$ with a variable number $t$ of types. Then the optimization problem
can be solved in time polynomial in $t$ and the binary encoding of $w^k,l^k,u^k,b^k,n_k$.
\et

The rest of the article is organized as follows. In Section 2 we discuss the
feasibility problem which is easier than the optimization problem and admits a more
efficient algorithm. In Section 3 we discuss the optimization problem, using the
results on feasibility. We conclude in Section 4 with further discussion of tables.

\section{Feasibility}

In this section we consider the feasibility problem for huge $n$-fold integer programs:
$$\mbox{is}\ \ \left\{x\in\Z^{dn}\ :\
A^{[n]}x=b\,,\ l\leq x\leq u\right\}\ \ \mbox{nonempty ?}$$

We begin with the case of {\em symmetric} programs, with one type, so that $t=1$,
over an $s\times d$ totally unimodular matrix $A$. So the data here consists of
the top right-hand side $a\in\Z^d$, and for all bricks the same lower and
upper bounds $l,u\in\Z_\infty^d$ and same right-hand side $b\in\Z^s$.
Then the set in question can be written as
\begin{equation}\label{symmetric_feasiblity}
\left\{x\in\Z^{dn}\ :\ \sum_{i=1}^nx^i=a
\,,\ Ax^i=b\,,\ l\leq x^i\leq u\,,\ i=1,\dots,n\right\}\ .
\end{equation}
We have the following lemma.

\bl{Symmetric_Lemma}
Let $A$ be totally unimodular. Then the set in (\ref{symmetric_feasiblity})
is nonempty if and only if $Aa=nb$ and $nl\leq a\leq nu$, and this can be
decided in time that is polynomial in the binary encoding of $n,l,u,a,b$,
even when $A$ is a variable part of the input.
\el
\bpr Suppose first that the set in (\ref{symmetric_feasiblity}) contains a feasible
point $x=(x^1,\dots,x^n)$. Then $Aa=A\sum_{i=1}^nx^i=\sum_{i=1}^nAx^i=nb$,
and $nl\leq a=\sum_{i=1}^nx^i\leq nu$. For the converse we use induction on $n$.
Suppose $a$ satisfies the conditions. If $n=1$ then $x^1:=a$ is a feasible point
in (\ref{symmetric_feasiblity}). Suppose now $n\geq 2$. Consider the system
$$l\leq y\leq u\,,\ \ Ay=b\,,\ \ (n-1)l\leq a-y\leq(n-1)u$$
in the variable vector $y$. Then $y={1\over n}a$ is a real solution to this system,
and therefore, since $A$ is totally unimodular, there is also an integer solution $x^n$
to this system. In particular, $Ax^n=b$ and $l\leq x^n\leq u$. Let ${\bar a}:=a-x^n$.
Then $A{\bar a}=A(a-x^n)=(n-1)b$ and $(n-1)l\leq {\bar a}=a-x^n\leq(n-1)u$.
It therefore now follows by induction that there is an integer solution
$(x^1,\dots,x^{n-1})$ to the $(n-1)$-fold program
$$\left\{x\in\Z^{d(n-1)}\ :\ \sum_{i=1}^{n-1}x^i={\bar a}
\,,\ Ax^i=b\,,\ l\leq x^i\leq u\,,\ i=1,\dots,n-1\right\}\ .$$
Then $\sum_{i=1}^nx^i={\bar a}+x^n=a$ and therefore $x:=(x^1,\dots,x^{n-1},x^n)$
is a feasible point in (\ref{symmetric_feasiblity}).
The statement about the computational complexity is obvious.
\epr

We proceed with the general case of $t$ types.
So the data now consists of $b^0\in\Z^d$ and for $k=1,\dots,t$,
lower and upper bounds $l^k,u^k\in\Z_\infty^d$, right-hand side $b^k\in\Z^s$,
and positive integer $n_k$, with $n_1+\cdots+n_t=n$. We denote by
$I_1\uplus\cdots\uplus I_t=\{1,\dots,n\}$ the natural partition
with $|I_k|=n_k$. So the set in question can be now written as
\begin{equation}\label{types_feasiblity}
\left\{x\in\Z^{dn}\ :\ \sum_{i=1}^nx^i=b^0
\,,\ Ax^i=b^k\,,\ l^k\leq x^i\leq u^k\,,\ k=1,\dots,t\,,\ i\in I_k\right\}\ .
\end{equation}

We have the following theorem asserting that when $A$ is totally unimodular the
feasibility problem is decidable in polynomial time even if the number $t$ of types
is variable. The algorithm underlying the proof uses only classical $n$-fold integer
programming and avoids the heavy results of \cite{GR} on integer cones used in \cite{Onn2}.

\bt{feasibility_theorem}
Let $A$ be a fixed totally unimodular matrix and consider the huge $n$-fold
program over $A$ with variable number $t$ of types.
Then it is decidable in time polynomial in $t$ and the binary encoding of
$l^k,u^k,b^k,n_k$, if the set in (\ref{types_feasiblity}) is nonempty.
\et
\bpr
Consider the following set of points of a classical $t$-fold integer program:
\begin{equation}\label{auxiliary_feasiblity}
\left\{y\in\Z^{dt}\ :\ \sum_{k=1}^ty^k=b^0
\,,\ Ay^k=n_kb^k\,,\ n_kl^k\leq y^k\leq n_ku^k\,,\ k=1,\dots,t\right\}\ .
\end{equation}
We claim that (\ref{types_feasiblity}) is nonempty if and only if
(\ref{auxiliary_feasiblity}) is nonempty, which can be decided within the claimed
time complexity by Proposition \ref{Classical} on classical $n$-fold theory.

So it remains to prove the claim. First, suppose $x$ is in (\ref{types_feasiblity}).
Define $y$ by setting $y^k:=\sum_{i\in I_k}x^i$ for $k=1,\dots,t$. Then we have
$\sum_{k=1}^ty^k=\sum_{i=1}^nx^i=b^0$, $Ay^k=A\sum_{i\in I_k}x^i=n_kb^k$,
and $n_kl^k\leq y^k=\sum_{i\in I_k}x^i\leq n_ku^k$, so $y$ is in
(\ref{auxiliary_feasiblity}). Conversely, suppose $y$ is in (\ref{auxiliary_feasiblity}).
For $k=1,\dots,t$ consider the symmetric $n_k$-fold program
$$\left\{(x^i:i\in I_k)\in\Z^{dn_k}\ :\ \sum_{i\in I_k}x^i=y^k
\,,\ Ax^i=b^k\,,\ l^k\leq x^i\leq u^k\,,\ i\in I_k\right\}\ .$$
Since $y$ is in (\ref{auxiliary_feasiblity}) we have that $Ay^k=n_kb^k$ and
$n_kl^k\leq y^k\leq n_ku^k$. Therefore, by Lemma \ref{Symmetric_Lemma}, this
program is feasible and has a solution $(x^i:i\in I_k)$.
Let $x=(x^1,\dots,x^n)$ be obtained by combining the solutions of these $t$ programs.
Then we have $\sum_{i=1}^nx^i=\sum_{k=1}^ty^k=b^0$ and
$Ax^i=b^k$ and $l^k\leq x^i\leq u^k$ for $k=1,\dots,t$ and $i\in I_k$, so $x$ is in
(\ref{types_feasiblity}). This completes the proof of the claim and the theorem.
\epr

\section{Optimization}

In this section we consider the optimization problem for huge $n$-fold programs:
$$\min\left\{\sum_{k=1}^t\sum_{i\in I_k}w^kx^i\, :\, x\in\Z^{dn},\sum_{i=1}^nx^i=b^0,
Ax^i=b^k,l^k\leq x^i\leq u^k,k=1,\dots,t,i\in I_k\right\}\ .$$
The optimization problem is harder than the feasibility problem in that we need to
actually produce an optimal solution if one exists. Since the problem is huge, meaning
that $n$ is encoded in binary, we cannot explicitly even write down a single point
$x\in\Z^{dn}$ in polynomial time. But it turns out that we can present $x$ compactly
as follows. For $k=1,\dots,t$ the set of all possible bricks of type $k$ is the following
$$S^k\ :=\ \{z\in\Z^d\ :\ Az=b^k\,,\ l^k\leq z\leq u^k\}\ .$$
We assume for simplicity that $S^k$ is finite for all $k$, which is the
case in most applications, such as in multiway table problems.
Let $\lambda^k:=(\lambda^k_z:z\in S^k)$ be a nonnegative integer tuple with entries
indexed by points of $S^k$. Each feasible point $x=(x^1,\dots,x^n)$ gives rise to
$\lambda^1,\dots,\lambda^t$ satisfying $\sum\{\lambda^k_z:z\in S^k\}=n_k$, where
$\lambda^k_z$ is the number of bricks of $x$ of type $k$ which are equal to $z$. Let
the {\em support} of $\lambda^k$ be $\supp(\lambda^k):=\{z\in S^k:\lambda^k_z\neq 0\}$.
Then a {\em compact presentation of $x$} consists of the restrictions of $\lambda^k$
to $\supp(\lambda^k)$ for all $k$. While the cardinality of $S^k$ may be exponential
in the binary encoding of the data $b^k,l^k,u^k$, it turns out that
a compact presentation of polynomial size always exists.
The following theorem was shown in \cite{Onn2} using the recent computationally heavy
algorithm of \cite{GR} which builds on \cite{ES}.

\bp{Fixed}
For fixed $d$ and $t$, the huge $n$-fold integer optimization problem with $t$ types,
over an $s\times d$ matrix $A$ which is part of the input, can be solved in
polynomial time. That is, in time polynomial in the binary encoding of
$A,l^k,u^k,b^k,n_k$, it can either be asserted that the problem is infeasible,
or a compact presentation $\lambda^1,\dots,\lambda^t$ of an optimal
solution with $|\supp(\lambda^k)|\leq 2^d$ for $k=1,\dots,t$ be computed.
\ep

We now show that for a totally unimodular matrix, we can solve the huge problem
even for variable $t$, extending both the above result and classical $n$-fold theory.

\vskip.2cm\noindent{\bf Theorem \ref{optimization_theorem}}
{\em Let $A$ be a fixed totally unimodular matrix and consider the huge \break $n$-fold
program over $A$ with a variable number $t$ of types. Then the optimization problem
can be solved in time polynomial in $t$ and the binary encoding of $w^k,l^k,u^k,b^k,n_k$.}

\vskip.2cm
\bpr
Consider the following classical $t$-fold integer optimization problem:
$$\min\left\{\sum_{k=1}^t w^ky^k\, :\, y\in\Z^{dt},\, \sum_{k=1}^ty^k=b^0,
\, Ay^k=n_kb^k,\, n_kl^k\leq y^k\leq n_ku^k,\, k=1,\dots,t\right\}.$$
By Proposition \ref{Classical} on classical $n$-fold theory we can either
assert the problem is infeasible, or obtain an optimal solution $y$, within the
claimed time complexity. As shown in the proof of Theorem \ref{feasibility_theorem},
if this problem is infeasible, then so is the original program, and we are done.
So assume we have obtained an optimal solution $y$.

For $k=1,\dots,t$ consider the symmetric $n_k$-fold program
$$\min\left\{\sum_{i\in I_k}w^kx^i\, :\, (x^i:i\in I_k)\in\Z^{dn_k}\,,\
\sum_{i\in I_k}x^i=y^k\,,\ Ax^i=b^k\,,\ l^k\leq x^i\leq u^k\,,\ i\in I_k\right\}\,.$$
As shown in the proof of Theorem \ref{feasibility_theorem}, this program is feasible.
Since this is a huge symmetric program, that is, with a single type, by
Proposition \ref{Fixed} we can compute in polynomial time a compact presentation
$\lambda^k$ with $|\supp(\lambda^k)|\leq 2^d$ of an optimal solution
$(x^i:i\in I_k)\in\Z^{dn_k}$. (In fact, any point in that program has the same objective
function value $\sum_{i\in I_k}w^kx^i=w^ky^k$ and is optimal to that program.)
Then $\lambda^1,\dots,\lambda^t$ obtained from all these programs provide
a compact presentation of a point $x=(x^1,\dots,x^n)$ feasible in the original program.
We claim this $x$ is optimal. Suppose indirectly there is a better point $\bar x$
and define $\bar y$ by ${\bar y}^k=\sum_{i\in I_k}{\bar x}^i$ for all $k$.
Then we have
$$\sum_{k=1}^t w^k{\bar y}^k\ =\ \sum_{k=1}^t\sum_{i\in I_k}w^k{\bar x}^i
<\sum_{k=1}^t\sum_{i\in I_k}w^kx^i=\sum_{k=1}^t w^ky^k\ ,$$
contradicting the optimality of $y$ in the $t$-fold program.
So indeed $\lambda^1,\dots,\lambda^t$ provide a compact presentation
of an optimal solution of the given huge $n$-fold program.
\epr

We make the following remark. The algorithm of Theorem \ref{feasibility_theorem}
for the feasibility problem involves only one application of the classical
$n$-fold integer programming algorithm of Proposition \ref{Classical}.
In contrast, the algorithm of Theorem \ref{optimization_theorem} is much heavier,
and in addition to one application of classical $n$-fold integer programming, uses $t$
times the algorithm of Proposition \ref{Fixed} for huge $n$-fold integer programming
with one type, which in turn uses the heavy algorithm for integer cones of \cite{GR}.

\section{Tables}

We now return to tables. Consider first $3$-way $l\times m\times n$ tables.
Index each table as $x=(x^1,\dots,x^n)$ with $x^i=(x^i_{1,1},\dots,x^i_{l,m})$.
Then the table problem is the $n$-fold program with matrix $A_{l,m}^{[n]}$
with $A_{l,m}$ the vertex-edge incidence matrix of the bipartite graph $K_{l,m}$.
Indeed, then $\sum_{i=1}^n x^i=g$ provides the vertical line-sum equations, and
$A_{l,m}x^i=b^k$ with $b^k=(e^k,f^k)$ provides the column and row sum equations
for $i\in I_k$. Since $A_{l,m}$ is totally unimodular,
Theorem \ref{optimization_theorem} implies our following claimed result.

\vskip.2cm\noindent{\bf Theorem \ref{ThreeWayLineSum}}
{\em The huge $3$-way table problem with a variable number $t$ of types can be solved in
time which is polynomial in $t$ and in the binary encoding of $w^k,e^k,f^k,g,n_k$.
In particular, deciding if there is a huge table with variable number of types is in P.}

\vskip.2cm
Let us continue with $4$-way $k\times l\times m\times n$ tables. Index each table as\break
$x=(x^1,\dots,x^n)$ with each $x^k$ an $k\times l\times m$ layer. Then  the table
problem is the $n$-fold program with matrix $A^{[n]}$ where $A=A_{k,l}^{[m]}$.
Now, unfortunately, for $k,l,m\geq 3$, the matrix $A$ is {\em not} totally unimodular.
Therefore, the results of the previous sections do not apply, and we remain with
the results of \cite{Onn2}, which are as follows.

\bp{FourWayLineSum}
The huge $4$-way table problem with fixed number $t$ of types is solvable in
time polynomial in the binary encoding of $w^k,n_k$ and the line sums. Moreover,
deciding if there is a huge table with $t$ variable is in NP intersect coNP.
\ep
The contrast between Theorem \ref{ThreeWayLineSum} and
Proposition \ref{FourWayLineSum} motivates the following.

\vskip.2cm\noindent
{\bf Open problem.}
What is the complexity of deciding feasibility of huge $4$-way tables
with a variable number of types ? In particular, is it in P for $3\times 3\times 3\times n$ tables ?

\section*{Acknowledgments}

The research of the first author was partially supported by a VPR Grant
at the Technion and by the Fund for the Promotion of Research at the Technion.
The research of the second author was partially supported by a
Bourse d'Aide \`a la Mobilit\'e Internationale from Universit\'e Paris Est.


\begin{thebibliography}{}

\bibitem{AT}
Aoki, S., Takemura, A.:
Minimal basis for connected Markov chain over $3\times3\times K$
contingency tables with fixed two-dimensional marginals.
Australian and New Zealand Journal of Statistics 45:229--249 (2003)

\bibitem{DHOW}
De Loera, J., Hemmecke, R., Onn, S., Weismantel, R.:
N-fold integer programming.
Discrete Optimization 5:231--241 (2008)

\bibitem{DO1}
De Loera, J., Onn, S.:
The complexity of three-way statistical tables.
SIAM Journal on Computing 33:819--836 (2004)

\bibitem{DO2}
De Loera, J., Onn, S.:
All linear and integer programs are slim 3-way transportation programs.
SIAM Journal on Optimization 17:806--821 (2006)

\bibitem{ES}
Eisenbrand, F., Shmonin, G.:
Carath\'eodory bounds for integer cones.
Operations Research Letters 34:564--568 (2006)

\bibitem{GR}
Goemans, M.X., Rothvo\ss, T.:
Polynomiality for Bin Packing with a Constant Number of Item Types.
Symposium on Discrete Algorithms 25:830--839 (2014)

\bibitem{HOR}
Hemmecke, R., Onn, S., Romanchuk, L.:
N-fold integer programming in cubic time.
Mathematical Programming 137:325--341 (2013)

\bibitem{HOW}
Hemmecke, R., Onn, S., Weismantel, R.:
A polynomial oracle-time algorithm for convex integer minimization.
Mathematical Programming 126:97--117 (2011)

\bibitem{HS}
Ho\c sten, S., Sullivant, S.:
Finiteness theorems for Markov bases of hierarchical models.
Journal of Combinatorial Theory Series A 114:311--321 (2007)

\bibitem{Onn}
Onn, S.: Nonlinear Discrete Optimization.
Zurich Lectures in Advanced Mathematics,
European Mathematical Society (2010),
available online at: {\tt http://ie.technion.ac.il/$\sim$onn/Book/NDO.pdf}

\bibitem{Onn2}
Onn, S.: Huge multiway table problems.
Discrete Optimization 14:72--77 (2014)

\bibitem{SS}
Santos, F., Sturmfels, B.:
Higher Lawrence configurations.
Journal of Combinatorial Theory Series A 103:151--164 (2003)

\end{thebibliography}
\end{document}